\documentclass[letterpaper,oneside,10pt]{article}

\usepackage[USenglish]{babel} 
\usepackage[T1]{fontenc}

\usepackage{lmodern} 

\usepackage{graphicx} 

\usepackage{float}
\usepackage{caption}
\usepackage{amsmath}
\usepackage{amsthm}

\usepackage{amsfonts}
\usepackage{amssymb}
\usepackage{mathtools}
\usepackage[hyphens]{url}

\let\Right\right 
\let\Left\left 
\makeatletter 
\def\right#1{\Right#1\@ifnextchar){\!\right}{}} 
\def\left#1{\Left#1\@ifnextchar({\!\left}{}} 
\makeatother

\begin{document}

\pagestyle{empty} 


\title{On the special harmonic numbers $H_{\lfloor p/9 \rfloor}$ and $H_{\lfloor p/18 \rfloor}$ modulo $p$}
\author{John Blythe Dobson (j.dobson@uwinnipeg.ca)}

\maketitle


\pagestyle{plain} 


\begin{abstract}
\noindent
Building on work of Zhi-Hong Sun, we establish congruences for the special harmonic numbers $H_\lfloor p/9 \rfloor$ and $H_{\lfloor p/18 \rfloor}$ modulo $p$, which contain respectively three and four distinct arithmetic components. We also obtain a complete determination modulo $p$ of the corresponding families of sums of reciprocals of the type studied by Dilcher and Skula. Applications to the first case of Fermat's Last Theorem are considered.

\end{abstract}

\section{Introduction}

\noindent
In recent years, considerable progress has been made in the determination modulo a prime $p$ of sums of the form

\begin{equation} \label{eq:SkulaSumDefinition}
s(k, N) :=
\sum_{j=\lfloor\frac{kp}{N}\rfloor + 1}^{\lfloor\frac{(k + 1)p}{N}\rfloor}
\frac{1}{j},
\end{equation}

\noindent
for integers $N$ and $k$ with $1 \le N < p$ and $0 \le k \le N-1$; the provision $j \neq p$ is necessary when $k=N-1$. In fact we need only consider cases with $k < N/2$ because by symmetry  $s(k, N) \equiv -s(N-1-k, N)$ mod $p$. These sums, probably first introduced by Lerch \cite{Lerch}, owe much of their interest to the 1995 result of Dilcher and Skula \cite{DilcherSkula1995} that for an exponent $p$ to be an exception to the first case of Fermat's Last Theorem (FLT), $p$ must divide every such sum $s(k, N)$ with $N \le 46$, and this interest has endured despite the publication of Wiles's complete proof of FLT the same year. In the case $N=0$, these sums can be conveniently written as a special kind of Harmonic number, $s(0, N) := H_{\lfloor p/N \rfloor}$.

The general problem of characterizing the sums in (\ref{eq:SkulaSumDefinition}) modulo $p$ was solved in two papers by the brothers Zhi-Hong Sun and Zhi-Wei Sun of 1992 and 2002 respectively (\cite{ZHSun1992}, \cite{ZWSun2002}), and an elementary exposition of their method in the case $k = 0$ has been given by the present writer \cite{DobsonRamus}. This approach does not explicitly reveal all the distinguishable components in the representations, as we propose to do here in the cases $N=9$ and $N=18$. Nonetheless, the 1992 paper (\cite{ZHSun1992}, pt. 2, Theorem 2.7) makes a crucial beginning in this direction with its evaluation of $s(1, 18)$, expressed in an alternate form. As this result has rarely if ever been cited, it is considered in some detail in the next section. From there we progress to the introduction of a companion result for $H_{\lfloor p/9 \rfloor}$, and finally to a complete set of calculations of $s(k, 9)$ and $s(k, 18)$ modulo $p$ for all values of $k$.

These final calculations are routine, and depend only on the well-known linear relations established either by Lerch \cite{Lerch} or by Dilcher and Skula \cite{DilcherSkula2011}, in which $q_p(b) := (b^{p-1} - 1)/p$ represents the Fermat quotient of $p$ in the base $b$:

\begin{equation} \label{eq:LerchDilcherSkula}
\begin{split}
s(0, 18) + s(2, 18) + s(4, 18) + s(6, 18) + s(8, 18) & \equiv -10 \cdot q_p(2) \pmod{p} \\
s(1, 18) + s(3, 18) + s(5, 18) + s(7, 18) & \equiv 8 \cdot q_p(2) \pmod{p} \\
s(0, 18) + 2 \cdot s(1, 18) + s(8, 18) & \equiv 0 \pmod{p}.
\end{split}
\end{equation}

\noindent
These are complemented by three relations drawn from the work of Emma Lehmer \cite{Lehmer1938}:

\begin{equation} \label{eq:Lehmer}
\begin{split}
s(0, 18) + s(1, 18) + s(2, 18) & = s(0, 6) \equiv -2 \cdot q_p(2) - \frac{3}{2} \cdot q_p(3) \pmod{p} \\
s(3, 18) + s(4, 18) + s(5, 18) & = s(1, 6) \equiv \phantom{-}2 \cdot q_p(2) \pmod{p} \\
s(6, 18) + s(7, 18) + s(8, 18) & = s(2, 6) \equiv -2 \cdot q_p(2) + \frac{3}{2} \cdot q_p(3) \pmod{p}.
\end{split}
\end{equation}

\noindent
Hendel \cite{Hendel} showed that $q_p(2)$ is an essential component of all $s(k, N)$ of even $N$, with the pattern observable in (\ref{eq:Lehmer}) persisting throughout.

\section{Toward the development of $H_{\lfloor p/18 \rfloor}$}

\noindent
Sun (\cite{ZHSun1992}, pt. 2, Theorem 2.7) begins with the development of a quotient which we shall designate $Z(p)$, such that, in effect,

\begin{equation} \label{eq:18Sun}
Z(p) := f(p)/p \equiv s(1, 18) - 2q_p(2) \pmod{p},
\end{equation}

\noindent
where $f(\cdot)$ is a sequence developed below. Where we write $s(1, 18)$, Sun (\cite{ZHSun1992}, pt. 2, Theorem 2.7) actually writes

\begin{equation*}
\sum_{k=1}^{\lfloor p/9 \rfloor} \frac{(-1)^{k-1}}{k},
\end{equation*}

\noindent
which is $-s^{\ast}(0, 9)$ in the notation of \cite{DobsonRamus}. Now if $T(N, m)$ represents the lacunary sum of binomial coefficients

\begin{equation*}
\sum_{\substack{j=0 \\ j \equiv 0 \bmod{N}}}^{m} \binom{m}{j},
\end{equation*}

\noindent
then the general result of Sun's three-part paper of 1992--1995 \cite{ZHSun1992}, with $f(\cdot)$ as defined in (\ref{eq:18Sun}), gives

\begin{equation} \label{eq:18formula}
f(p) \equiv 9 \cdot T(9, p) - 2 \cdot T(2, p) - 7 \pmod{p^2}.
\end{equation}

\noindent
This implies that $f(\cdot)$ satisfies the ninth-order linear recurrence (not actually written out by Sun),

\begin{equation} \label{eq:Z0}
\begin{split}
f(n) = & 8f(n - 1) - 29f(n - 2) + 62f(n - 3) - 86f(n - 4) + 80(n - 5) \\
       & - 50f(n - 6) + 20f(n - 7) - 5f(n - 8) + f(n - 9).
\end{split}
\end{equation}

\noindent
The first 36 terms of this sequence are listed in the table below.

\begin{table} [H]
\begin{center}
\begin{tabular}{ r | r || r | r || r | r }
$n$ & $f(n)$ & $n$ & $f(n)$ & $n$ & $f(n)$ \\
\hline
 1 &       0 & 13 &    $-1755$ & 25 &   $-10841355$ \\
 2 &    $-2$ & 14 &       1636 & 26 &   $-24927437$ \\
 3 &    $-6$ & 15 &      12279 & 27 &   $-49854867$ \\
 4 &   $-14$ & 16 &      37426 & 28 &   $-88165112$ \\
 5 &   $-30$ & 17 &      87720 & 29 &  $-135356601$ \\
 6 &   $-62$ & 18 &     175447 & 30 &  $-166501907$ \\
 7 &  $-126$ & 19 &     307287 & 31 &  $-109482111$ \\
 8 &  $-254$ & 20 &     464776 & 32 &     202204690 \\
 9 &  $-501$ & 21 &     560190 & 33 &    1101562302 \\
10 &  $-932$ & 22 &     348313 & 34 &    3176261536 \\
11 & $-1551$ & 23 &  $-731055$ & 35 &    7325660004 \\
12 & $-2114$ & 24 & $-3798314$ & 36 &   14651320015 \\
\end{tabular}
\caption{Values of $f(n)$ for $n \le 36$}
\label{Table_1}
\end{center}
\end{table}

\noindent
Sun establishes that this sequence satisfies $f(p) \equiv 0$ mod $p$, so that $f(p) \equiv f(1) \equiv 0$ and $f(\cdot)$ is both a Fermat sequence \cite{Gillespie} and a prime-divisible sequence \cite{Minton}. Now if $f(\cdot)$ were subjected to multisection in the usual manner, being split according to the residue class of $n$ modulo $9$, the multisections would be characterized by the fifth-order recurrence

\begin{equation}
f(n) = -246f(n-1) + 13606f(n-2) + 245f(n-3) - 13605f(n-4) + f(n-5).
\end{equation}

\noindent
But Sun does not take this standard approach. Rather, he first notices that the even terms of $f(\cdot)$ can be increased by 2 and the terms with $n \equiv 3 \bmod{6}$ increased by 6 without affecting the terms where $n$ is a prime $> 3$, yielding a new sequence that can be divided throughout by 3, finally giving

\begin{equation} \label{eq:g}
g(n) = 0, 0, 0, -4, -10, -20, -42, -84, -165, -310, -517, -704, \dots,
\end{equation}

\noindent
which is represented by the ninth-order recurrence

\begin{equation} \label{eq:Z2}
\begin{split}
g(n) = & 6g(n - 1) - 15g(n - 2) + 20g(n - 3) - 18g(n - 4) + 18(n - 5) \\
       & - 20g(n - 6) + 12g(n - 7) - 3g(n - 8) + g(n - 9).
\end{split}
\end{equation}

\noindent
It may seem counterproductive to transform $f(\cdot)$ into another sequence of the same order, but $g(\cdot)$ possesses a very interesting property. The sequence is amenable to an unusual kind of multisection, in which the selected terms in each subsequence come not from a single residue class of $n$ modulo $p$, but rather from two complementary residue classes of $n$:

\begin{equation*} \label{eq:ZCasesStartingValues}
\begin{dcases}
0, -84, -310, 29240, 102429, -8309145, -29388370, \dots & \text{for $n \equiv \pm 1 \pmod{9}$} \\
0, -42, -517, 12476, 154926, -3613785, -45118867, \dots & \text{for $n \equiv \pm 2 \pmod{9}$} \\
0, -20, -704, 4095, 186732, -1266104, -55500635,  \dots & \text{for $n \equiv \pm 3 \pmod{9}$} \\
-4, -10, -585, 546, 116105, -243685, -36494037,   \dots & \text{for $n \equiv \pm 4 \pmod{9}$}.\\
\end{dcases}
\end{equation*}

\noindent
(Although the case of $n \equiv \pm 3 \bmod{9}$ obviously includes no prime $> 3$, we include it for the sake of completeness; however the simpler sequence formed by the terms with $n$ divisible by 9 is not important for our purposes.)

With the above sets of seven starting values, the four subsequences all satisfy the same seventh-order linear recurrence

\begin{equation} \label{eq:Z}
\begin{split}
h(m) = & h(m - 1) - 246h(m - 2) + 246h(m - 3) + 13605h(m - 4) \\
       & - 13605h(m - 5) - h(m - 6) + h(m - 7),
\end{split}
\end{equation}

\noindent
where

\begin{equation*} \label{eq:ZCasesSeed}
m = 
\begin{dcases}
2 \left\lfloor \frac{n}{9} \right\rfloor + 1 & \text{if $n \equiv 1, 2, 4 \pmod{9}$} \\
2 \left\lfloor \frac{n}{9} \right\rfloor + 2 & \text{if $n \equiv -1, -2, -4 \pmod{9}$}. \\
\end{dcases}
\end{equation*}

\noindent
Sun next devises a decomposition of $g(\cdot)$ such that

\begin{equation} \label{eq:SunCases}
g(n) = \begin{dcases}
2j(n) - 1                        & \text{if $n \equiv \phantom{\pm} 0 \pmod{9}$} \\
j(n + 1) - 2                     & \text{if $n \equiv \pm 1 \pmod{9}$} \\
j(n) - j(n - 1) - 2              & \text{if $n \equiv \pm 2 \pmod{9}$} \\
-j(n) - 1                        & \text{if $n \equiv \pm 3 \pmod{9}$} \\
-j(n + 1) - j(n) + j(n - 1) - 2  & \text{if $n \equiv \pm 4 \pmod{9}$}, \\
\end{dcases}
\end{equation}

\noindent
and finds that $j(\cdot)$ satisfies a lacunary recurrence of the third order,

\begin{equation} \label{eq:SunSeries}
\begin{split}
j(1) & = 0, j(2) = 2, j (3) = -1 \\
j(n) & = 3j(n-2) - j(n-3) \quad (n \ge 3).
\end{split}
\end{equation}

\noindent
The sequence $j(\cdot)$ retains the original property of being a Fermat sequence and a prime-divisible, since $j(p) \equiv j(1) \equiv 0$ mod $p$. So finally, considering only the cases pertinent to primes $p > 9$, we have

\begin{equation} \label{eq:SunFinalDefinition}
Z(p) := 3 \cdot \frac{g(p)}{p} \equiv \begin{dcases}
3 \cdot \frac{j(p+1) - 2}{p}                  & \text{if $p \equiv \pm 1 \pmod{9}$} \\
3 \cdot \frac{j(p) - j(p-1) - 2}{p}           & \text{if $p \equiv \pm 2 \pmod{9}$} \\
3 \cdot \frac{-j(p+1) - j(p) + j(p-1) - 2}{p} & \text{if $p \equiv \pm 4 \pmod{9}$}. \\
\end{dcases}
\end{equation}

\section{The development of $H_{\lfloor p/9 \rfloor}$}

\noindent
Prompted by the work of Hendel \cite{Hendel}, we now introduce a companion sequence to Sun's $Z(p)$. The applicability of this sequence to the present problem will be justified in detail in the next section, but first we prefer to dispose of the technical details of its construction. We begin with the quotient

\begin{equation} \label{eq:9standard}
X(p) := k(p)/p \equiv s(0, 9) - s(0, 3) \pmod{p}.
\end{equation}

\noindent
Following the method of the Sun brothers, this satisfies

\begin{equation} \label{eq:9formula}
k(p) \equiv 3 \left\{-6 \cdot T(18, p) + 3 \cdot T(9, p) + 2 \cdot T(6, p) - T(3, p) + 2 \right\} \pmod{p^2},
\end{equation}

\noindent
where $T(N, m)$ is again the lacunary sum of binomial coefficients

\begin{equation*}
\sum_{\substack{j=0 \\ j \equiv 0 \bmod{N}}}^{m} \binom{m}{j}.
\end{equation*}

\noindent
The divisibility of $k(p)$ by $p$ follows from the work of Sun, and since $k(p) \equiv k(1) \equiv 0 \bmod{p}$, $k(\cdot)$ is both a Fermat sequence and a prime-divisible sequence. If we consider $k(n)/3$ for all natural $n$, we have a sequence of which the first 36 terms are listed in the table below.

\begin{table} [H]
\begin{center}
\begin{tabular}{ r | r || r | r || r | r }
$n$ & $k(n)/3$ & $n$ & $k(n)/3$ & $n$ & $k(n)/3$ \\
\hline
 1 &     0 & 13 &    2873 & 25 &        5218265 \\
 2 &     0 & 14 &    6734 & 26 &        5218265 \\
 3 &  $-1$ & 15 &   15014 & 27 &              2 \\
 4 &  $-4$ & 16 &   32132 & 28 &    $-20242870$ \\
 5 & $-10$ & 17 &   66368 & 29 &    $-78528607$ \\
 6 & $-19$ & 18 &  132734 & 30 &   $-226111984$ \\
 7 & $-28$ & 19 &  257393 & 31 &   $-572533513$ \\
 8 & $-28$ & 20 &  483626 & 32 &  $-1343905180$ \\
 9 &     2 & 21 &  877799 & 33 &  $-2992450957$ \\
10 &   110 & 22 & 1529363 & 34 &  $-6395344954$ \\
11 &   407 & 23 & 2527769 & 35 & $-13201132948$ \\
12 &  1145 & 24 & 3873017 & 36 & $-26402265898$ \\
\end{tabular}
\caption{Values of $k(n)/3$ for $n \le 36$}
\label{Table_2}
\end{center}
\end{table}

\noindent
The oscillating sequence $k(n)$ returns to 2, the value of the constant term in (\ref{eq:9formula}), whenever $n$ is an odd multiple of 9; and apart from the first two terms, the $n$th term is positive when the residue of $n \bmod{36}$ lies between 9 and 27, and negative when the residue lies between 0 and 8 or between 28 and 35. The sequence is defined by the seven initial values stated in the table above and a linear recurrence of the seventh order,

\begin{equation} \label{eq:9recurrence}
\begin{split}
k(n) = & 7k(n - 1) - 21k(n - 2) + 36k(n - 3) - 39k(n - 4) \\
       & + 27k(n - 5) - 12k(n - 6) + 3k(n - 7).
\end{split}
\end{equation}

\noindent
A multisection $r(n)$ of $k(n)$ based on any residue class of $n$ modulo 18 satisfies the fourth-order recurrence

\begin{equation} \label{eq:9recurrenceMultisectioned}
r(n) = -199097r(n-1) - 18108279r(n-2) + 18287694r(n-3) - 19683r(n-4),
\end{equation}

\noindent
which implies that (\ref{eq:9recurrence}) can be reformulated as the lacunary recurrence

\begin{equation} \label{eq:9recurrenceLacunary}
\begin{split}
k(n) = & -199097k(n - 18) - 18108279k(n - 36) \\
       & + 18287694k(n - 54) + 19683k(n - 72).
\end{split}
\end{equation}

\section{Expressions for $s(k, 9)$ and $s(k, 18)$}

\noindent
From the definition (\ref{eq:9standard}) and the classic result of Glaisher (\cite{Glaisher}, p.\ 50) that $s(0, 3) \equiv -\frac{3}{2}q_p(3) \bmod{p}$, we have

\begin{equation} \label{eq:9definition}
s(0, 9) \equiv X(p) - \frac{3}{2}q_p(3) \pmod{p}.
\end{equation}

\noindent
Comparing this with Sun's evaluation of $s(1, 18)$, we immediately infer that $s(0, 18) \equiv -2q_p(2) - \frac{3}{2}q_p(3) + X(p) - Z(p)$. This and the remaining results below satisfy the linear relations gathered in the Introduction. In view of the relations (\ref{eq:18Sun}) and (\ref{eq:9definition}), as an added check all representations were verified numerically for known $p$ for which $q_p(2)$ vanishes modulo $p$ (the Wieferich primes, OEIS A001220), for which $q_p(3)$ vanishes modulo $p$ (the Mirimanoff primes, OEIS A014127), and for which $X(p)$ or Z(p) vanishes modulo $p$. This establishes that the components are truly independent, and that any other valid definitions of $X(p)$ and $Z(p)$ could differ from ours only by a scaling factor.

The full sets of values follow:

\begin{table} [H]
\begin{center}
\begin{tabular}{@{}r | rrrr@{}}
$k$ & $q_p(3)$ & $X(p)$ & $Z(p)$ \\
\hline
0 &  $-\frac{3}{2}$ &    1 &   0  \\
1 &            $-3$ & $-3$ & $-1$ \\
2 &               3 &    2 &   1  \\
3 &               0 &    1 &   1  \\
4 &               0 &    0 &   0  \\
\end{tabular}
\caption{Coefficients of the components of $s(k, 9)$}
\label{Table_3}
\end{center}
\end{table}

\begin{table} [H]
\begin{center}
\begin{tabular}{@{}r | rrrr@{}}
$k$ & $q_p(2)$ & $q_p(3)$ & $X(p)$ & $Z(p)$ \\
\hline
0 & $-2$ &  $-\frac{3}{2}$ &    1 & $-1$ \\
1 &    2 &               0 &    0 &    1 \\
2 & $-2$ &               0 & $-1$ &    0 \\
3 &    2 &            $-3$ & $-2$ & $-1$ \\
4 & $-2$ &               6 &    6 &    4 \\
5 &    2 &            $-3$ & $-4$ & $-3$ \\
6 & $-2$ &            $-6$ & $-5$ & $-2$ \\
7 &    2 &               6 &    6 &    3 \\
8 & $-2$ &   $\frac{3}{2}$ & $-1$ & $-1$ \\
\end{tabular}
\caption{Coefficients of the components of $s(k, 18)$}
\label{Table_4}
\end{center}
\end{table}

\noindent
In light of the proof by Dilcher and Skula \cite{DilcherSkula1995} that a prime $p$ failing the first case of FLT must vanish simultaneously for $s(k, 9)$ and $s(k, 18)$ with all values of $k$, the results of a search for vanishing values of these sums may be of interest:

\begin{table} [H]
\begin{center}
\begin{tabular}{ r | l }
$k$ & $p$ \\
\hline
0  & 677, 6691, 532199813         \\
1  & 151, 457, 971, 1439, 12613   \\
2  & 241, 739, 37799, 3112456729  \\
3  & 97, 58193                    \\
\end{tabular}
\caption{Values of $p < 5,330,000,000,000$ for which $s(k, 9)$ vanishes modulo $p$}
\label{Table_5}
\end{center}
\end{table}

\begin{table} [H]
\begin{center}
\begin{tabular}{ r | l }
$k$ & $p$ \\
\hline
0  & 5235774727037                                \\
1  & 47, 1777, 217337                             \\
2  & 167                                          \\
3  & 1171, 37783, 28525219                        \\
4  & 137, 251, 1087, 1301, 2111, 5749, 428687393  \\
5  & 4177, 1581479                                \\
6  & 108541, 48303223                             \\
7  & 149, 35267                                   \\
8  & ---                                          \\
\end{tabular}
\caption{Values of $p < 5,330,000,000,000$ for which $s(k, 18)$ vanishes modulo $p$}
\label{Table_6}
\end{center}
\end{table}

\noindent
These calculations were performed in PARI/GP using standard matrix-powerng methods to evaluate the various sequences modulo $p$, and $s(8, 18)$ is the only remaining case for which no result has been found. It will be noted that no value of $p$ appears in more than one position.

\section{Connection with the first case of Fermat's Last Theorem}

\noindent
For a prime $p$ to fail the first case of Fermat's Last Theorem (FLT), Wieferich proved that $q_p(2)$ must vanish modulo $p$, while Mirimanoff proved that $q_p(3)$ must vanish modulo $p$. Thus, the relations (\ref{eq:18Sun}) and (\ref{eq:9definition}) and the 1995 result of Dilcher and Skula \cite{DilcherSkula1995} mentioned in the Introduction imply the vanishing of both $Z(p)$ and $X(p)$ modulo $p$ as necessary conditions for the failure of the first case of FLT. A computer search to $p < 5,330,000,000,000$ found only the following cases where the expressions vanish modulo $p$:

For $Z(p)$: $179$, $1949$, $28885849$.

For $X(p)$: $2$, $13$, $19$, $2423$.

\noindent
Within the range of $p$ tested, there is no $p$ for which $s(k, 18)$ vanishes simultaneously for distinct values of $k$, nor any value of $k$ for which any two of the three quantities $s(k, 18)$, $Z(p)$, and $X(p)$ vanish simultaneously.

\section{Conclusion}

\noindent
The 1995 paper of Dilcher and Skula \cite{DilcherSkula1995} considers $s(k, N)$ with $N$ ranging from $2$ through $46$. Computational studies of these numbers have focused mainly on the case $k=0$, or $H_{\lfloor p/N \rfloor}$. Although Dilcher and Skula state (pp.\ 389--390) that $s(0, N)$ vanishes modulo $p$ for some value of $p < 2000$ for every $N$ in this range except $N = 5$, in fully $22$ of these cases ($N = 7, 11, 12, 13, 14, 15, 16, 17, 18, 20, 22, 28, 29, 30, 31, 35, 38, 39, 43, 44, 45, 46$) this statement is true only if the sums are permitted to be vacuous, which surely cannot have been the authors' intention. Both in previous work \cite{DobsonSpecialHarmonicNumbers} and in the present paper we have disregarded vacuous sums and extended the search in those cases to higher $p$. For $N$ between $2$ and $46$, at least one solution $p$ is now known for all but the cases $N = 5, 12, 20, 29, 31, 43$.

\end{document}